\documentclass{amsart}
\usepackage{color}
\newtheorem{Theorem}{Theorem}[section]
\newtheorem{Definition}[Theorem]{Definition}
\newtheorem{Proposition}[Theorem]{Proposition}

\newtheorem{Lemma}[Theorem]{Lemma}
\newtheorem{Corollary}[Theorem]{Corollary}
\theoremstyle{remark}

\newtheorem{Example}[Theorem]{Example}

\newtheorem{Assertion}[Theorem]{Assertion}

\def\il{\int\limits_}

\def\eps{\varepsilon}

\def\ovr{\overline}

\def\om{\omega}
\def\Om{\Omega}
\def\al{\alpha}
\def\gm{\gamma}
\def\th{\theta}

\def\dl{\delta}

\def\bd{\partial}
\def\lm{\lambda}

\def\sm{\setminus}
\def\sbs{\subset}

\def\re{{\mathbf {Re\,}}}

\def\be{\begin{enumerate}}
\def\ee{\end{enumerate}}
\def\bT{\begin{Theorem}}
\def\eT{\end{Theorem}}
\def\bP{\begin{Proposition}}
\def\eP{\end{Proposition}}
\def\bD{\begin{Definition}}
\def\eD{\end{Definition}}
\def\bE{\begin{Example}}
\def\eE{\end{Example}}
\def\bL{\begin{Lemma}}
\def\eL{\end{Lemma}}
\def\bC{\begin{Corollary}}
\def\eC{\end{Corollary}}
\def\bA{\begin{Assertion}}
\def\eA{\end{Assertion}}
\def\A{{\mathcal A}}

\def\oD{\ovr{\mathbb D}}
\def\aD{\mathbb D}
\def\aN{\mathbb N}

\def\aC{\mathbb C}

\begin{document}

\title{ Plurisubharmonically separable complex manifolds}
\author{Evgeny A. Poletsky and Nikolay Shcherbina}
 \begin{abstract} Let $M$ be a complex manifold and $PSH^{cb}(M)$ be the
 space of bounded continuous plurisubharmonic functions on $M$. In this
 paper we study when the functions from $PSH^{cb}(M)$ separate points.
 Our
 main results show that this property is equivalent to each of the
 following properties of $M$:\be
\item the core of $M$ is empty.
\item for every $w_0\in M$ there is a continuous plurisubharmonic
    function $u$ with the logarithmic singularity at $w_0$.
\ee
\par Moreover, the core of $M$ is the disjoint union of the sets $E_j$
that are 1-pseudoconcave in the sense of Rothstein and have the
following
Liouville property: every function from $PSH^{cb}(M)$ is constant on
each
of $E_j$.
\end{abstract}
\thanks{The first author was partially supported by a grant from Simons
Foundation.}
\keywords{plurisubharmonic functions, 1-pseudoconcave sets.}
\subjclass[2010]{Primary 32U05; Secondary 32F10, 32U35.}
\address{Department of Mathematics,  Syracuse University, 215 Carnegie
Hall, Syracuse, NY 13244, USA}
\email{eapolets@syr.edu}
\address{Department of Mathematics, University of Wuppertal, 42119
Wuppertal, Germany}
\email{shcherbina@math.uni-wuppertal.de}
\maketitle
\section{Introduction}
\par Let $M$ be a complex manifold and $PSH^{cb}(M)$ be the space of
bounded continuous plurisubharmonic functions on $M$. In this paper we
study obstructions to separation of points  in a complex manifold $M$ by
functions from $PSH^{cb}(M)$.
\par There are complex manifolds, for example compact manifolds, where all
plurisubharmonic functions are constants. There are parabolic manifolds,
for example $\aC^n$, where all bounded plurisubharmonic functions are
constants. For their characterization see \cite{R2} and \cite{AS}. But
there are plenty of complex manifolds like $\aD\times\aC$, $\aD$ is the
unit disk, where the space $PSH^{cb}(M)$ is large but, nevertheless, does
not separate points.
\par The first main results of our paper can be summarized in the
following
theorem.

\noindent

{\bf Theorem I.} {\it For a complex manifold $M$ the following statements
are equivalent: \be
\item the functions from $PSH^{cb}(M)$ separate points of $M$;
\item for every point $w_0\in M$ there is a function $u\in PSH^{cb}(M)$
    that is smooth and strictly plurisubharmonic near $w_0$ ;
\item for every point $w_0\in M$ there are a negative continuous
    plurisubharmonic function $v$ on $M$ and constants $C_1$ and $C_2$
    such
    that $\log|z-w_0|+C_1<v(z)<\log|z-w_0|+C_2$ near $w_0$.
\ee
}

\par The main obstruction to separation is the  set $\mathbf{c}(M)$ of all
points $w\in M$, where every function of $PSH^{cb}(M)$ fails to be smooth
and strictly plurisubharmonic near $w$. It was first introduced and
systematically studied by Harz--Shcherbina--Tomassini in
\cite{HST1}--\cite{HST3} and was called the {\it core} of $M$. Observe
that
directly from the definition one concludes that $\mathbf{c}(M)$ is a
closed
subset of $M$.  Among the main properties of the core established in these
papers we mention here the following result   that will be one of the
important technical tools in the present paper.

\noindent
{\bf Theorem.} (see \cite[Theorem 3.2]{HST1}) {\it Let $M$ be a complex
manifold. Then the set $\mathbf{c}(M)$ is $1$-pseudo\-concave in the sense
of Rothstein. In particular, $\mathbf{c}(M)$ is pseudoconcave in $M$ if
$\dim_{\mathbb C} M = 2$.}

\par Our next main result is the following theorem that was proved in
\cite{HST2} when the dimension of $M$ is two.

\noindent
{\bf Theorem II.} {\it Let $M$ be a complex manifold. Then the set
$\mathbf{c}(M)$ is the disjoint union of the sets $E_j$, $j\in J$, that
are
1-pseudoconcave in the sense of Rothstein and have the following Liouville
property: every function from $PSH^{cb}(M)$ is constant on each of $E_j$.}

\par The equivalence of statements 2 and 3 in Theorem I is proved in
Theorem \ref{T:3.2}. The equivalence of 1 and 2 and Theorem II are proved
in  Theorem \ref{T:4.6}. Note that the sets $E_j$ can have a fractal
nature
(see \cite{HST4}).

\par Observe that if $E \subset M$ is a set that is 1-pseudoconcave in the
sense of Rothstein in $M$ which has the property that each bounded above
continuous plurisubharmonic function $\varphi$ on $M$ is constant on
$E$, then $E \subset \mathbf{c}(M)$ (for details see Lemma 3.1 in
\cite{HST2}). Hence, mentioned above Theorem II clarify the phenomenon of
both, existence and the structure of the core.

\par Note that in this paper we are dealing mainly with the core defined
using continuous plurisubharmonic on $M$ functions (the core
$\mathbf{c}^0(M)$ in terminology of \cite{HST3}), while the main object of
the study in \cite{HST1}-\cite{HST3} was the core $\mathbf{c}(M)$ defined
using smooth plurisubharmonic on $M$ functions. Note also, that another
proof of  Theorem II for cores defined by smooth plurisubharmonic
functions
was obtained  by Slodkowski \cite{S2} using essentially different methods.
His proof also covers the case of {\it minimal kernels} which are defined
and studied in \cite{ST}.

\section{Basic facts}

\par Let $PSH(M)$ be the space of plurisubharmonic functions on a complex
manifold $M$. We denote by $PSH^b(M)$ ($PSH^c(M)$) the subspaces of
bounded
above (continuous) functions in $PSH(M)$. Let $PSH^{cb}(M)$ be the
intersection of $PSH^b(M)$  and $PSH^c(M)$.

\par The {\it plurisubharmonic envelope} $u$ of a function $\psi$ on a
complex manifold $M$ is the supremum of all plurisubharmonic functions on
$M$ that are less or equal to $\psi$ on $M$. If $\psi$ is upper
semicontinuous, then by the Poletsky--Rosay theorem (see \cite{P} and
\cite{R1})
\[u(z)=\inf_{f\in\A(\oD,M),f(0)=z}\frac1{2\pi}\il0^{2\pi}\psi(f(e^{i\th}))\,d\th.\]

\par Given a complex manifold $M$, an open set $U\sbs M$ and a point
$w_0\in M$ let $v(z)$ be the supremum of all negative plurisubharmonic
functions on $M$ that are less than $-1$ on $U$. The {\it harmonic measure
} $\om(w_0,U,M)$ of the set $U$ at $w_0$ is equal to $-v(w_0)$. By the
Poletsky--Rosay theorem $\om(w_0,U,M)$ is the supremum over all analytic
disks $f:\,\oD\to M$, $f(0)=w_0$, of the normalized measure $\lm$ of the
set $\{\th:\,0\le\th<2\pi, f(e^{i\th})\in U\}$.

\par A closed set $E\sbs M$ is called {\it 1-pseudoconcave in the sense of
Rothstein} if for any $z_0\in E$ and for any strictly plurisubharmonic
function $\rho$ defined on a neighborhood $V$ of $z_0$ at any neighborhood
$U\sbs\sbs V$ containing $z_0$ there is a point $z\in E\cap U$ where
$\rho(z)>\rho(z_0)$.

\par By their definition  these sets are perfect, i.e., have no isolated
points. But it may happen that such a set is compact. For example, take
the
unit ball $B$ in $\aC^2$ and blow-up a complex projective line $E$ at the
origin. We get a complex manifold $M$. Clearly, $E$ is a set that is
1-pseudoconcave in the sense of Rothstein.

\par But if $M$ is Stein, then any connected component $X$ of $E$ is
non-compact. Indeed, by \cite[Lemma 5.4]{HST2}, $X$ is a 1-pseudoconcave
in
the sense of Rothstein set. If $X$ is compact and $\rho$ is a smooth
strictly plurisubharmonic exhaustion function on $M$, then we take
the minimal $a$ such that $X\sbs\{\rho(z)\le a\}$. If $z_0\in
X\cap\{\rho=a\}$, then $\rho(z)\le\rho(z_0)$ for all $z\in X$ and it
shows that $X$ is not a set that is 1-pseudoconcave in the sense of
Rothstein.

\section{Cores and functions with logarithmic singularities}

\par Let $M$ be a connected complex manifold and $w_0\in M$. Let us choose
some local coordinates near $w_0$ and define the pluricomplex Green
function with a pole at $w_0$ as
\[g_M(z,w_0)=\sup\{u\in PSH(M):\,u<0 \text{ on }M, u(z)<\log|z-w_0|+c
\text{ near } w_0\}.\]

The function $g_M$ is a plurisubharmonic in $z$ on $M$ and
$g_M(w_0,w_0)=-\infty$.

\par Let $G_M(w_0,c)=\{z\in M:\,g(z,w_0)<c\}$, $c<0$, be Green balls.
Evidently  Green balls are open. They are also connected. Indeed, if the
open set $G_M(w_0,c)$ has a connected component $U$ that does not contain
$w_0$, then we take the function $u(z)$ equal to $c$ on $U$ and to
$g_M(z,w_0)$ on $M\sm U$. This function is, evidently, plurisubharmonic on
$M\sm\bd U$. But if $z\in\bd U$, then $g_M(z,w_0)\ge c$. Otherwise a
neighborhood of $z$ will be in $U$. Hence $u$ is upper semicontinuous and,
since $u\ge g_M$, plurisubharmonic. But $g_M(z,w_0)\ge u(z)$ on $M$ and,
consequently, $g_M(z,w_0)\ge c$ on $U$. Hence $U$ does not lie in
$G_M(w_0,c)$.

\par We will need the following lemma.

\noindent
\bL\label{L:3.1}  Let $u$ be a bounded plurisubharmonic function on $M$
that is smooth and strictly plurisubharmonic near $w_0\in M$. Then for any
coordinate neighborhood $V$ of $w_0$ there is a ball $U\sbs V$ centered at
$w_0$ and a pluriharmonic function $h$ on $U$ such that $u-h>0$ near $\bd
U$ and $u(w_0)-h(w_0)<0$.\eL

\begin{proof} Let us choose some coordinate neighborhood $X\sbs V$ around
$w_0$ such that $u$ is smooth on $X$. For $z\in X$ we set
\[L_{w_0}(z)=\sum_{i,j=1}^nu_{z_i,z_j}(w_0)(z_i-(w_0)_i)(z_j-(w_0)_j)\]
and
\[H_{w_0}(z)=\sum_{i,j=1}^nu_{z_i,\ovr z_j}(w_0)(z_i-(w_0)_i)(\ovr
z_j-(\ovr w_0)_j).\]
The Taylor expansion of $u$ at $w_0$ is
\[u(z)=u(w_0)+2\re (\nabla u(w_0),z-w_0) +\re
L_{w_0}(z)+\frac12H_{w_0}(z)+o(\|z-w_0\|^2).\]
\par The function
\[v(z)=u(w_0)+2\re (\nabla u(w_0),z-w_0) +\re L_{w_0}(z)\] is
pluriharmonic
on $X$ and there is a neighborhood $Y\sbs\sbs X$ such that the function
$u(z)-v(z)>c|z-w_0|^2$ on $Y$ for some $c>0$. We take as $U$ a ball of
radius $r$, centered at $w_0$ and compactly belonging to $Y$. Let
$a=cr^2/2$. If $h=v+a$, then the function $u-h>0$ on $\bd U$ and
$u(w_0)-h(w_0)<0$.
\end{proof}
\par Now we can prove the equivalence of statements 2 and 3 in Theorem I.

\bT\label{T:3.2} A point $w_0\not\in \mathbf{c}(M)$ if and only if  there
are a negative continuous plurisubharmonic function $v$ on $M$ and
constants $C_1$ and $C_2$ such that $\log|z-w_0|+C_1<v(z)<\log|z-w_0|+C_2$
near $w_0$.\eT
\begin{proof} Let $u$ be a continuous negative plurisubharmonic function
on
$M$ strictly plurisubharmonic near $w_0$. By Lemma \ref{L:3.1} in some
coordinate neighborhood $V$ of $w_0$ there is a ball $U$ of radius $r$
centered at $w_0$ and a pluriharmonic function $h$ on  $V$ such that
$u-h>0$ on $\bd U$ and $a=u(w_0)-h(w_0)<0$. There is another ball
$B\sbs\sbs U$ of radius $s$ centered at $w_0$ such that $u-h<a/2$ on $B$.
The pluricomplex Green function $g_U(z,w_0)=\log(|z-w_0|/r)$ on $U$.  We
take a constant $c>0$ such that $c\log(s/r)>a/2$. Then $cg_U>u-h$ on $\bd
B$. Define the function $v$ on $M$ as $g_U+h/c$ on $B$,
$\max\{g_U+h/c,u/c\}$ on $U\sm B$ and $u/c$ outside of $U$. Since
$g_U+h/c>u/c$ on $\bd B$ and $g_U+h/c<u/c$ on $\bd U$, the function $v$ is
continuous and plurisubharmonic on $M$.

\par To prove the converse statement we take in some coordinate
neighborhood of $w_0$ an open ball $B$ of radius $r$ centered at $w_0$
such that
$\log|z-w_0|+C_1<v(z)<\log|z-w_0|+C_2$
on $\ovr B$. Let $B_1=B(w_0,r/2)$. Then we take the function
$\phi(z)=a(|z-w_0|^2-r^2)$ on $B$, where $a>0$ is chosen so that
$-3ar^2/4<\log(r/2)+C_1$. This implies that $\phi<v$ on $\bd B_1$. We
define the function $u$ on $M$ as $\max\{v,\phi\}$ on $B_1$ and as $v$
everywhere else.  Clearly $u$ is a negative continuous plurisubharmonic
function on $M$.

\par Since $v(z)<\log|z-w_0|+C_2$ near $w_0$, there is an open
neighborhood
$U$ of $w_0$ where $v<\phi$. Hence $u=\phi$ on $U$ and is smooth and
strictly plurisubharmonic near $w_0$.
\end{proof}

\bC\label{C:3.3} If $w_0\not\in \mathbf{c}(M)$, then
$g_M(w,w_0)>\log|z-w_0|+C(w_0)$ in some
coordinate neighborhood of $w_0$ for some constant $C(w_0)$.\eC

\par  The following lemma tells us that the notions of cores and maximal
functions are related to each other.

\bL\label{L:3.4} Let $\Om\sbs\sbs M$ be a domain biholomorphic to a
bounded
domain in $\aC^n$, $w_0\in\Om$  and let $\phi\in PSH^{cb}(M)$. Let $u$ be
a
continuous plurisubharmonic function on $\Om$ such that $\limsup_{z\to
z_0,z_0\in\bd\Om}u(z)\le \phi(z_0)$. If $w_0\in \mathbf{c}(M)$, then
$u(w_0)\le\phi(w_0)$.\eL

\begin{proof} Suppose that $u(w_0)=\phi(w_0)+a$, $a>0$. Let $u_1=u-a/2$
and
let $U=\{z : u_1(z)>\phi(z)\}$. Then $U\sbs\sbs\Om$. We take an
open set $V\sbs\sbs\Om$ such that $U\sbs\sbs V$ and $u_1<\phi-b$ on
$\bd V$ for some $b>0$. Since $\Om$ is
biholomorphic to a domain in $\aC^n$, there is a decreasing sequence
of smooth plurisubharmonic functions $v_j$ on a neighborhood of $\ovr V$
converging to $u_1$ on $\ovr V$. If $j$ is sufficiently large, then we may
assume that $v_j(w_0)>\phi(w_0)+a/2$, but $v_j<\phi-b/2$ on $\bd V$.

\par The domain $\Om$ is biholomorphic to a bounded domain in $\aC^n$.
Hence there is a bounded, smooth and strictly  plurisubharmonic function
$v$ on $\Om$. We set $u_2=v_j+cv$, where $c>0$ is chosen in such a way
that
$u_2<\phi$ on $\bd V$ and $u_2(w_0)>\phi(w_0)$. We define the function
$\psi$ as $\max\{u_2,\phi\}$ on $V$ and as $\phi$ on $M\sm\ovr V$. Since
the function $\psi$ is continuous and plurisubharmonic on $M$, smooth and
strictly plurisubharmonic near $w_0$, we see that $w_0\not\in
\mathbf{c}(M)$.
\end{proof}

\bC\label{C:3.5}  If $\mathbf{c}(M)=M$, then any $\phi\in PSH^{cb}(M)$ is
maximal.\eC

\section{Foliation of cores}

\par Let $w_0$ be a point in $M$. A point $z\in A^b(w_0)$ ($A^{cb}(w_0)$)
if $u(z)\le u(w_0)$ for any $u\in PSH^{b}(M)$ ($u\in PSH^{cb}(M)$).
Evidently, $A^b(w_0)\sbs A^{cb}(w_0)$.

\par Let us list some easily derived properties of the sets $A^b(w_0)$ and
$A^{cb}(w_0)$.

\noindent
\bP\label{P:4.1} \be \item The set $A^{cb}(w_0)$ is closed.
\item If $z_0\in A^b(w_0)$ ($z_0\in A^{cb}(w_0)$), then $A^b(z_0)\sbs
    A^b(w_0)$ ($A^{cb}(z_0)\sbs A^{cb}(w_0)$).
\item A point $w\in A^b(w_0)$ if and only if $\om(w,U,M)=1$ for any
    neighborhood $U$ of $w_0$.
\item If $w\in A^b(w_0)$, then $g_M(w,w_0)=-\infty$. Hence the set
    $A^b(w_0)$ is pluripolar if $g_M(\cdot,w_0)\not\equiv-\infty$.
\item A point $w_0\not\in \mathbf{c}(M)$ if and only if
    $A^{cb}(w_0)=\{w_0\}$.
\item Any bounded holomorphic function on $M$ is constant on $A^b(w_0)$.
\ee\eP

\begin{proof} (1). If $u\in PSH^{cb}(M)$, then the set $\{w:\,u(w)\le
u(w_0)\}$ is closed. Since $A^{cb}(w_0)$ is the intersection of such sets
over all $u\in PSH^{cb}(M)$, (1) follows.

\par (2) is evident. To show (3) we note that if $w\in A^b(w_0)$ and $U$
is
a neighborhood of $w_0$, then any negative plurisubharmonic function on
$M$
that is less or equal to $-1$ on $U$ is also less or equal to $-1$ at $w$.
Hence $\om(w,U,M)=1$.
\par If $\om(w,U,M)=1$ for any neighborhood $U$ of $w_0$ and $u$ is a
plurisubharmonic function on $M$ such that $u\le A$ on $M$, then for any
$\eps>0$ we take a neighborhood $U$ of $w_0$, where $u\le u(w_0)+\eps$ on
$U$. In view of the mentioned above Poletsky-Rosay theorem, it follows
that there is an analytic disk $f$ such that $f(0)=w$ and the measure of
the set $E=\{\th:\,0\le\th<2\pi, f(e^{i\th})\in U\}$ is greater than
$1-\eps$. Hence $u(w)\le (u(w_0)+\eps)(1-\eps)+A\eps$. Since $\eps>0$ is
arbitrary, we see that $u(w)\le u(w_0)$.

\par Let us show (4). Since $g_M(w_0,w_0)=-\infty$, for any $N<0$ there is
a neighborhood $U$ of $w_0$ such that $g_M(z,w_0)\le N$ on $U$. If $w\in
A^b(w_0)$, then $\om(w,U,M)=1$ and for any $\eps>0$ there is an analytic
disk $f$ such that $f(0)=w$ and the measure of the set
$E=\{\th:\,0\le\th<2\pi, f(e^{i\th})\in U\}$ is greater than $1-\eps$.
Hence $g_M(f(0),w_0)\le (1-\eps)N$. Since $N$ and $\eps$ are arbitrary, we
see that $g_M(w,w_0)=-\infty$. Thus $g_M(w,w_0)=-\infty$ on $A^b(w_0)$ and
if $g_M(\cdot,w_0)\not\equiv-\infty$, then $A^b(w_0)$ is pluripolar.

\par (5). If a point $w_0\not\in \mathbf{c}(M)$, then by Theorem
\ref{T:3.2} there is a negative continuous function $u$ on $M$ such that
$\{u=-\infty\}=\{w_0\}$. Hence $A^{cb}(w_0)=\{w_0\}$.

\par If $A^{cb}(w_0)=\{w_0\}$, then we take a ball $B$ centered at
$w_0$ and for each point $z\in\bd B$ find a function $u_z\in PSH^{cb}(M)$
such that $u_z(w_0)=0$ and $u_z(z)>0$. We take finitely many points
$z_1,\dots,z_k$ in $\bd B$ such that the function
$u=\max\{u_{z_1},\dots,u_{z_k}\}$ is greater than $0$ on $\bd B$. Then the
maximal plurisubharmonic function $v$ on $B$ with boundary values equal to
$u$ has the value at $w_0$ strictly bigger than $u(w_0)$. By Lemma
\ref{L:3.4} $w_0\not\in \mathbf{c}(M)$.

\par (6). Let $f$ be a bounded holomorphic function on $M$ and $|f|<r-1$
on
$M$. If $w\in A^b(w_0)$, then $\log|f(w)+a|\le\log|f(w_0)+a|$ for any
$a\in\aC$ with $|a|=r$. If $h(w)=(f(w)+a)^{-1}$, then again
$\log|h(w)|\le\log|h(w_0)|$. Hence $|f(w)+re^{i\th}|=|f(w_0)+re^{i\th}|$
for any $\th\in[0,2\pi]$. Therefore
\[|f(w)|^2-|f(w_0)|^2=2\re((f(w_0)-f(w))re^{-i\th})\] and this implies
that
$f(w)=f(w_0)$.
\end{proof}
\par In Theorem \ref{T:cbr} below we prove the major property of the sets
$A^{cb}$. For its proof we need the following lemma.
\bL\label{L:nabla} Suppose that $A$ is a closed set in the closure $\ovr
B$ of the unit ball $B\sbs\aC^n$ centered at the origin, $0\in A$, and
$\rho$ is a strictly plurisubharmonic function defined on an open
neighborhood of $\ovr B$ such that $\rho(0)=0$ and $\rho\le 0$ on $A$. Then
there are $z_1\in A\cap B$ and a strictly plurisubharmonic function
$\rho_1$ defined on an open neighborhood $U$ of $z_1$  such that
$\rho_1(z_1)=0$, $\rho_1\le 0$ on $A\cap U$ and $\nabla\rho_1(z_1) \neq
0$.
\eL
\begin{proof} Since the function $\rho$ is strictly plurisubharmonic,
there is $\eps>0$ such that the function $\rho^*(z)= \rho(z) -\eps\|z\|^2$
is also strictly plurisubharmonic on $\ovr B$.
Observe that the function $\rho^*\le0$ on $A$, $\rho^*(0) = 0$ and,
moreover, now we also have that $\rho^* \le -\eps<0$ on $A\cap\bd B$.
If $\rho^*\le0$ on some open ball centered at the origin, then it is zero
on this ball and this is impossible since $\rho^*$ is strictly
plurisubharmonic. Hence we can take a point $a$ with $\rho^*(a) > 0$ so
close
to the origin that for all $t \in [0,1]$ we have $\rho*<-\eps/2<0$ on
$(A+ta)\cap\bd B$.
Consider now the function $F : [0,1] \rightarrow
[0,\infty)$ defined as
\[F(t)=\max\{\rho^*(z) : z \in(A+ta) \cap \ovr B\}.\]
It is obvious that the function $F$ is continuous and $[0,\rho^\ast(a)]
\subset F([0,1])$. By Sard's theorem there is a non-critical value
$s\in[0,\rho^*(a)]$ for the function $\rho^*$. Let $t^* \in [0,1]$ be
the value of the parameter $t$ such that $F(t^\ast) = s$ and let $z^*$ be a
point in $A+t^*a$, where $\rho^*$ attains its maximum equal to $s$. Since
$\rho^*<0$ on $(A+ta) \cap\bd B$, $z^*\in B$, and since $\rho^*(z^*)>0$ the
point $z_1=z^*- t^*a$ is also in $B$. Hence $z_1$ and the function
$\rho_1(z)=\rho^*(z+t^*a)-s$, satisfy all requirements of the lemma.
\end{proof}
\par Now we can prove the theorem.
\bT\label{T:cbr} If $A^{cb}(w_0)\ne\{w_0\}$, then the set $A^{cb}(w_0)$ is
1-pseudoconcave in the sense of Rothstein.\eT
\begin{proof} Suppose to get a contradiction that there are $z_0\in
A^{cb}(w_0)$, a neighborhood $U\sbs\sbs M$ of $z_0$ and a strictly
plurisubharmonic function $\rho$ defined on an open neighborhood of $\ovr
U$ such that the set $A^{cb}(w_0)\cap \ovr U\cap\{\rho>\rho(z_0)\}$ is
empty. By Lemma \ref{L:nabla} we may assume that in some local coordinates
$(z_1,\dots,z_n)$ the
set $U$ is a ball $B(0,1)$ in $\aC^n$, $z_0=0$, $\rho(0)=0$ and
$\nabla\rho(0)\ne0$. Then there is a biholomorphic change of coordinates 
near $z_0$ such that the function $\rho$ becomes strictly convex and $z_0$ 
is still equal to $0$. We may assume that $\nabla\rho(0)=(1,0,\dots,0)$.
\par First, we consider the case when $z_0\ne w_0$ and $w_0\not\in B$. Let
$h(z)=\re z_1$ and $S_0=\{\rho=0\}$.  There is a smaller ball $B=B(0,r)$
and $\dl>0$ such that $\rho(z+te_1)>0$ when $0<t\le\dl$ and $z\in S_0\cap
B(0,r)$. The Taylor expansion of $\rho$ at $z_0$ is
\[\rho(z)=2h(z)+Q(z)+o(\|z\|^2),\] where $Q(z)\ge a\|z\|^2$, $a>0$, is a
real quadratic form. If $\rho(z)<0$, then $a\|z\|^2\le
-2h(z)+o(\|z\|^2)$ and  there is $\eps>0$ such that the set $G=\{z\in
B:\,-\eps<h(z)<0,\rho(z)<0\}$ compactly belongs to $B$. We take
$0<\dl_1<\min\{\dl,\eps/2\}$ so small that the set
$G_1=G+(\dl_1,0,\dots,0)$ also compactly belongs to $B$. Let $S=\bd
G_1\cap\{S_0+\dl_1e_1\}$.
\par Observe that $\rho>0$ on $S$ and, therefore, $S\cap
A^{cb}(w_0)=\emptyset$. So for every $z\in S$ there is a function $u_z\in
PSH^{cb}(M)$ such that $u_z(w_0)=0$ and $u_z(z)>0$. Replacing $u_z$ with
$\max\{u_z,0\}$ we may assume that $u_z\ge 0$ on $M$. By compactness of
$S$, there are points $z_1,\dots,z_k$ on $S$ such that the function
$u=\max\{u_{z_1},\dots,u_{z_k}\}>\al>0$ on $S$.

\par Recall that $h=\dl_1-\eps<0$ on $\bd G_1\sm S$. If $\dl_2=\dl_1-\eps$
and $m$ is the infimum of $u$ on $S$, then the function
\[v(z)=\frac{m}{2\eps}(h(z)-\dl_2)\] is pluriharmonic on $\ovr
G_1$, $v=0$ on $\bd G_1\sm S$, $v>0$ on $G_1$ and
\[v(z)\le v(\dl_1,0,\dots,0)=m/2\] on $S$. Thus $v\le u$ on $\bd G_1$.

\par We define the function $u_1$ on $M$ as $\max\{u,v\}$ on $\ovr G_1$
and
$u$ on $M\sm\ovr G_1$. This function will be plurisubharmonic on $M$ and
$u_1(z_0)>0$. Hence $z_0\not\in A^{cb}(w_0)$ and we came to a
contradiction. Therefore, the set $A^{cb}(w_0)$ is 1-pseudoconcave in the
sense of Rothstein at any point $z_0\in A^{cb}(w_0)$, $z_0\ne w_0$.

\par If $z_0=w_0$ and there are a neighborhood $U$ of $z_0$ and a strictly
plurisubharmonic function $\rho$ defined on $U$ and such that the set
$A^{cb}(w_0)\cap U\cap\{\rho>\rho(w_0)\}$ is empty, then we perform the
same construction and note that $G_1$ is biholomorphic to a bounded domain
in $\aC^n$ containing $w_0$ and a plurisubharmonic function $v$ on $\ovr
G_1$ has its boundary values less or equal to $u$ while $v(w_0)>u(w_0)$.
Hence  by Lemma \ref{L:3.4}, $w_0\not\in\mathbf{c}(M)$ and by
(5) $A^{cb}(w_0)=\{w_0\}$.
\end{proof}

\par Let $A^{cb}_e(w_0)$ be the set of $w\in M$ such that $u(w)=u(w_0)$
for
any $u\in PSH^{cb}(M)$.

\bL\label{L:4.2} If $w_0\in \mathbf{c}(M)$, then there is a function $u\in
PSH^{cb}(M)$ equal to 0 on $A^{cb}_e(w_0)$ and negative on $A^{cb}(w_0)\sm
A^{cb}_e(w_0)$.\eL

\begin{proof} For any point $w\in A^{cb}(w_0)\sm A^{cb}_e(w_0)$  there is
a
function $u_w\in PSH^{cb}(M)$ such that $u_w(w_0)=0$ and $u_w(w)<0$.
Replacing $u_w$ with $\max\{u_w,-1\}$ we may assume that the uniform norm
$\|u_w\|$ of each function $u_w$ is finite.

\par Let us take a countable family of points $w_k\in A^b(w_0)\sm
A^{cb}_e(w_0)$ and their neighborhoods $U_k$ such that $\cup
U_k=A^{cb}(w_0)\sm A^{cb}_e(w_0)$ and $u_{w_k}<0$ on $U_k$. Note that
$u_{w_k}\le0$ on $A^{cb}(w_0)$ and $u_{w_k}=0$ on $A^{cb}_e(w_0)$.

\par Let $b_k=\|u_{w_k}\|$. We find $\gm_k>0$ such that the series
$\sum_{k=1}^\infty\gm_kb_k$ converges. Let
$u=\sum_{k=1}^\infty\gm_ku_{w_k}$. Then $u\in PSH^{cb}(M)$, $u=0$ on
$A^{cb}_e(w_0)$ and $u<0$ on $A^{cb}(w_0)\sm A^{cb}_e(w_0)$.
\end{proof}

\par The following lemma is the analog of \cite[Lemma 2.2]{S1} that can be
used on complex manifolds and is proved as in the original.

\bL\label{L:4.3} \it Let $K$ be a compact set in a complex manifold $M$.
Suppose that there is a bounded smooth strictly  plurisubharmonic function
$\rho$ defined on a neighborhood $V$ of $K$. Let $u$ be a plurisubharmonic
function on $V$. Assume that there is a non-empty compact set $L\sbs K$
such that $\max_Ku>\max_Lu$. Then there are a point $w^*$ in $K\sm L$, a
neighborhood $U$ of $w^*$ and a smooth strictly plurisubharmonic function
$v$ on
$U$ such that $v(w^*)=0$ and $v(w)<0$ when $w\ne w^*$ belongs to $K\cap
U$.\eL

\begin{proof} There is $\eps>0$ such that $M=\max_Ku_\eps>\max_Lu_\eps$
for
the function $u_\eps=u+2\eps\rho$. Take $w_1\in K$ such that
$u_\eps(w_1)=M$. Then $w_1\in K\sm L$.

\par  We can use the same argument as in Lemma \ref{L:3.1} and find a
sufficiently small coordinate ball $U$ centered at $w_1$ and a
pluriharmonic
function $h$ on $U$ such that the function $\rho_1=\rho+h$ attains its
strict
minimum at $w_1$ and $\rho_1(w_1)=0$. Let $v_0=u_\eps-\eps\rho_1-M$. Then
$v_0(w_1)=0$, $v_0(w)<0$  when $w\ne w_1$ belongs to $K\cap U$, and
$v_0<a<0$ on $\bd U$. Let us take a decreasing sequence of smooth
plurisubharmonic functions $u_j$, $j\in\aN$, converging to $u$ on an open
neighborhood of $\ovr U$. Then the strictly plurisubharmonic functions
$v_j=u_j+\eps(\rho+h)-M$ decreasingly converge to $v_0$. Hence there is $j$
such that $v_j<a/2$ on $\bd U$ while $v_j(w_1)\ge 0$. We take $w^*\in K\cap
U$, where $v_j$ attains its maximum $M'$ on $K\cap\ovr U$ and let
$v(w)=v_j(w)-M'-\dl\|w-w^*\|^2$, where $\dl>0$ is a sufficiently small
number.
\end{proof}

\par Following the terminology of \cite{S1} we say that a closed set $X$
in
a complex manifold $M$ has the {\it local maximum property} if $X$ is
perfect (i.e. has no isolated points) and for any $z_0\in X$ there is
an open neighborhood  $V$ of $z_0$ in $M$ with compact closure such that
if
an open set $U\sbs\sbs V$ contains $z_0$ and the set $L=X\cap\bd U$ is
non-empty, then
\[\sup_{X\cap\ovr U}u\le\sup_L u\]
for any  plurisubharmonic function $u$ on $V$.

\bT\label{T:4.4} If $M$ is a complex manifold, then any closed set $X\sbs
M$ has the local maximum property if and only if it is $1$-pseudoconcave
in
the sense of Rothstein.\eT

\begin{proof} Let $X$ be a set that is $1$-pseudoconcave in the sense of
Rothstein and $w_0\in X$. Take an open coordinate neighborhood $V$ of
$w_0$. Suppose that $U\sbs\sbs V$ is an open set containing $w_0$ and the
set $L=X\cap\bd U$ is non-empty. If  $u$ is a plurisubharmonic function on
$V$ such that $\sup_{X\cap\ovr U}u>\sup_L u$, then by Lemma \ref{L:4.3}
there are a point $w^*$ in $X\cap U$, a neighborhood $W$ of $w^*$ and a
strictly plurisubharmonic function $v$ on $W$ such that $v(w^*)=0$ and
$v(w)<0$ when $w\ne w^*$ belongs to $X\cap W$. But this contradicts the
assumption that the set $X$ is a $1$-pseudoconcave in the sense of
Rothstein. Hence $X$ has the local maximum property.

\par Conversely, if $X$ has the local maximum property and is not
$1$-pseudoconcave in the sense of Rothstein, then there is a point $w_0\in
X$ and a strictly plurisubharmonic function $\rho$ defined on a
neighborhood $V$ of $w_0$ such that $\rho(w)\le\rho(w_0)$ when $w\in X\cap
V$. We may assume that $V$ lies in a coordinate ball centered at $w_0$ and
find $\eps>0$ such the function $u(w)=\rho(w)-\eps|w-w_0|^2$ is still
strictly plurisubharmonic. Since $X$ is perfect there is an open set
$U\sbs\sbs V$ that contains $w_0$ such that the set $L=X\cap\bd U$ is
non-empty. Hence, $\sup_Lu<u(z_0)$ and this contradiction shows that $X$
is
$1$-pseudoconcave in the sense of Rothstein.
\end{proof}

\bT\label{T:4.5} If $M$ is a complex manifold and $w_0\in \mathbf{c}(M)$,
then the set $A^{cb}_e(w_0)$ has the local maximum property.\eT

\begin{proof} It is clear that the set $A^{cb}_e(w_0)$ is closed. Since
$A^{cb}_e(w_0)=A^{cb}_e(w_1)$ for any $w_1\in A^{cb}_e(w_0)$, it
suffice to prove the theorem at $w_0$. Let $U$ be an open set with compact
closure containing $w_0$ and the set $L=A^{cb}(w_0)\cap\bd U$ is
non-empty.
Let us show that the set $A^{cb}_e(w_0)\cap\bd U$ is also non-empty.
\par Since, in view of Theorem \ref{T:cbr}, the set $A^{cb}(w_0)$ is
$1$-pseudoconcave in the sense of Rothstein, it has the local maxi\-mum
property. By Lemma \ref{L:4.2} there is a continuous function $u$ on $M$
equal to 0 on $A^{cb}_e(w_0)$ and negative on $A^{cb}(w_0)\sm
A^{cb}_e(w_0)$. If the set $A^{cb}_e(w_0)\cap\bd U$ is empty, then
$\sup_{A^{cb}(w_0)\cap\ovr U}u>\sup_Lu$ and this contradicts to the local
maximum property of $A^{cb}(w_0)$. Hence $A^{cb}_e(w_0)\cap\bd
U\ne\emptyset$ and since the set  $A^{cb}(w_0)$ is perfect, the set
$A^{cb}_e(w_0)$ is also perfect.
\par Now let $V\sbs M$ be an open set with compact closure, an open set
$U\sbs\sbs V$ contains $w_0$ and the set $A^{cb}_e(w_0)\cap\bd U$ is
non-empty. Suppose that there is a plurisubharmonic function $v$ defined
on
$V$ such that $\sup_{A^{cb}_e(w_0)\cap\ovr U}v>\sup_{A^{cb}_e(w_0)\cap\bd
U} v$. We may assume that $\sup_{A^{cb}_e(w_0)\cap\ovr V}v=0$ so the
function $v$ is negative on $A^{cb}_e(w_0)\cap\bd U$. There is a
neighborhood $W$ of the latter set such that $v$ is negative on this
neighborhood. Let $u$ be a function from Lemma \ref{L:4.2}. The set
$K=(A^{cb}(w_0)\cap\bd U)\sm W$ is compact and does not contain points of
$A^{cb}_e(w_0)$ and, therefore, there is a constant $a>0$ such that the
function $v_1=v+au$ is negative on this set. Since $v$ is negative on $W$
and $u$ is non-positive on $A^{cb}(w_0)$, the function $v_1\le-\eps$
on $A^{cb}(w_0)\cap\bd U$ for some $\eps>0$. Also
$\sup_{A^{cb}_e(w_0)\cap\ovr U}v_1=0$ because $u=0$ on $A^{cb}_e(w_0)$.
Hence $\sup_{A^{cb}(w_0)\cap\ovr U}v_1>\sup_{A^{cb}(w_0)\cap\bd U}v_1$.
But
the set $A^{cb}(w_0)$ has the local maximum property and this means that
we
got a contradiction. Therefore, the set $A^{cb}_e(w_0)$ has the local
maximum property.
\end{proof}

\par  As a corollary to the latter result we obtain the following theorem.

\noindent
\bT\label{T:4.6} Let $M$ be a complex manifold with non-empty core
$\mathbf{c}(M)$. Then:\be
\item for every $w_0\in M$ the set $A^{cb}_e(w_0)\ne\{w_0\}$ if and only
    if
    $w_0\in\mathbf{c}(M)$;
\item if $w_0\in\mathbf{c}(M)$, then the set $A^{cb}_e(w_0)$ is a is
    1-pseudoconcave in the sense of Rothstein, lies in $\mathbf{c}(M)$ and
    all functions $u \in PSH^{cb}(M)$ are constants on $A^{cb}_e(w_0)$;
\item the core $\mathbf{c}(M)$ of $M$ can be decomposed into the disjoint
    union of closed sets $E_j$, $j\in J$, that are 1-pseudoconcave in the
    sense of Rothstein and have the following Liouville property: Every
    function $\varphi
\in PSH^{cb}(M)$ is constant on each of the sets $E_j$.
\ee
\eT

\begin{proof} (1) If $w_0\in \mathbf{c}(M)$, then by Theorem
\ref{T:4.5} the set $A^{cb}_e(w_0)$ has the local maximum property and,
therefore, is not equal to $\{w_0\}$. Conversely, if $w_0\not\in
\mathbf{c}(M)$, then $A^{cb}(w_0)=\{w_0\}$ by Proposition
\ref{P:4.1}(5) and, consequently, $A^{cb}_e(w_0)=\{w_0\}$.

\par (2) By Theorems \ref{T:4.4} and \ref{T:4.5} the set $A^{cb}_e(w_0)$
is
$1$-pseudoconcave in the sense of Rothstein. Let $w_1\in A^{cb}_e(w_0)$
and
$w_1\ne w_0$. Then $w_0\in A^{cb}_e(w_1)$. Hence $A^{cb}_e(w_1)\ne\{w_1\}$
and $w_1\in \mathbf{c}(M)$. By the definition of the sets $A^{cb}_e(M)$
all
functions $u \in PSH^{cb}(M)$ are constants on $A^{cb}_e(w_0)$.

\par (3) If $A^{cb}_e(w_0)\cap A^{cb}_e(w_1)\ne\emptyset$, then
$A^{cb}_e(w_0)=A^{cb}_e(w_1)$. Hence the relation $w_0\sim w_1$ if $w_1\in
A^{cb}_e(w_0)$ is an equivalence relation. Let $J=\mathbf{c}(M)/\sim$
and $[w]_\sim$ be the equivalence class of $w\in\mathbf{c}(M)$. Then
$\mathbf{c}(M)$ is the disjoint union of the sets $E_j=A^{cb}_e(w)$,
$[w]_\sim=j$.
\end{proof}

\par If the manifold $M$ is not Stein, it may happen that the sets
$A^{cb}(w_0)$ and $A^{cb}_e(w_0)$ are compact. For example, take the unit
ball $B$ in $\aC^2$ and blow-up a complex projective line $X$ at the
origin. We get a complex manifold $M$ and a holomorphic mapping $F$ of $M$
onto $B$ such that $F(X)=\{0\}$. If $u$ is a plurisubharmonic function on
$B$, then the function $v=u\circ F$ is plurisubharmonic on $M$ while
any plurisubharmonic function on $M$ is constant on $X$. Hence
$\mathbf{c}(M)=X$ and $A^{cb}(w_0)=A^{cb}_e(w_0)=X$ for $w_0 \in X$.

\par The following theorem shows that such phenomenon can exist only on
non-Stein manifolds.

\bT\label{T:4.7} Let $M$ be a complex manifold and let $X\ne\{w_0\}$ be
the
closed connected component of $A^{cb}(w_0)$ containing $w_0$. If $X$ is
compact, then $A^{cb}(w_0)=X$ and $M$ is not Stein.\eT

\begin{proof} If $X$ is compact, we can take a decreasing sequence of open
neighbourhoods $V_p$ of $X$ with compact closure such that
$\bigcap_{p=1}^{\infty}V_p = X$ and $\bd V_p \cap A^{cb}(w_0) = \emptyset$
for every $p \in \mathbb{N}$ (for the existence of the sets $V_p$ see, for
example, Proposition 2 on page 108 of [NN]). If we fix a set $V_p$ from
this sequence, then, since points $w\in\bd V_p$ do not belong to
$A^{cb}(w_0)$, for each of them there is a function $u_w\in PSH^{cb}(M)$
such that $u_w(w)>0$ and $u_w(w_0)=0$. Note that $u_w\le0$ on
$A^{cb}(w_0)$. Let us pick up finitely many points $w_k$, $1\le k\le m$,
on $\bd V_p$ such that the function $u=\max\{u_{w_1},\dots,u_{w_m}\}$ is
positive on $\bd V_p$. It follows from the compactness of $\bd V_p$
that $u>\al$ on $\bd V_p$ for some $\al>0$.

\par Let us take the function $v$ on $M$ that is equal to $u$ on $V_p$ and
to $\max\{u,\al/2\}$ on $M\sm V_p$. This function will be plurisubharmonic
and $v(w)>v(w_0)=0$ for $w\in M\sm V_p$. Hence $A^{cb}(w_0)\sbs V_p$. Then
finally from $\bigcap_{p=1}^{\infty}V_p = X$ we conclude that
$A^{cb}(w_0)=X$.
\end{proof}

\section{Example}
\par Here we present an example which refutes the following possible
conjectures about the sets $\mathbf{c}(M), A^b(w_0)$ and $A^{cb}(w_0)$:
\be\item the set $A^b(w_0)$ is closed;
\item the set $A^b(w_0)=A^{cb}(w_0)$;
\item if $w_1\in A^b(w_0)$, then $w_0\in A^b(w_1)$;
\item if $w_0\in\mathbf{c}(M)$, then there is a set $E$ containing
    $w_0$ and at least one more point such that any $u\in PSH^b(M)$ is
    constant on $E$.
\ee
\par Let
\[\phi(\zeta)=\sum_{k=1}^\infty\gm_k\log|(\zeta-k^{-1})/2|,\] where all
$\gm_k>0$ are chosen so that the series converges for all $\zeta\in\aC$,
$\zeta\ne k^{-1}$, $k\in\aN$, $\phi(0)=-3$ and $\phi(\zeta)<-2$ on the
unit
disk $\aD$. Clearly the function $\phi$ is subharmonic. Consider the
domain
\[\Om=\{(z_1,z_2)\in\aC^2:\,\psi(z_1,z_2)<0\},\] where
\[\psi(z_1,z_2)=|z_1|^2+|z_2|^2+\log|(z_1-1/2)/2|+\phi(z_2).\]
This domain contains the bidisk $\aD^2$ and the lines
$L_k=\{z_2=k^{-1}\}$,
$k\in\aN$, and the line $L_0=\{z_1=1/2\}$. If $u$ is a bounded above
plurisubharmonic function on $\Om$, then $u$ is constant on each of these
lines. Hence by connectedness, it is constant on the union $X$ of this
lines.

\par Moreover, if $X_0=\{(z_1,0)\in\Om\}$, i.e., if
$|z_1|^2+\log|(z_1-1/2)/2|<3$, then $e^\psi$ is continuous and bounded
on $\Om$ and it is smooth and strictly plurisubharmonic outside of
$X_0\cup X$. This shows that $\mathbf{c}(\Om)\sbs X_0\cup X$. On the
other hand, clearly, $X\sbs \mathbf{c}(\Om)$ and, since cores are
closed, it follows that $X_0\cup X\sbs \mathbf{c}(\Om)$.
\par To refute the conjecture (4) we take $w_0=(0,0)$. If
$w\not\in\mathbf{c}(\Om)$, then, by Corollary \ref{C:3.3}, the
function $g_M(w,w)=-\infty$ while $g_M(w_0,w)>-\infty$. If $w\in X$, then
$\psi(w)=-\infty$ while $\psi(w_0)=-3-\log 4$. If $w\in X_0$, then the
family
\[\psi_{a, b}(z_1,z_2)=a|z_1 - b|^2+|z_2|^2+\log|(z_1-1/2)/2|+\phi(z_2)\]
for  properly chosen $0<a<1$ and $b \in \mathbb{C}$ will  consist of the
bounded above on $\Om$ plurisubharmonic functions that separate the points
of $X_0$ from $w_0$.

\par If $u\in PSH^b(\Om)$, then for any $\eps>0$ there is a
neighborhood $U$ of $(0,0)$ such that $u(z_1,z_2)<u(0,0)+\eps$ on $U$.
Since $U$ meets some lines $z_2=k^{-1}$, we see that $u<u(0,0)+\eps$
on $X$. Thus $u\le u(0,0)$ on $X$ by arbitrariness of $\eps>0$ and,
thereofore, $X\sbs A^b((0,0))$. However, for the function $\psi$ we have
$\psi < 0$ on $\Om$, $\psi(0, 0) = - 3 - \log4$ and $\psi(z_1, k^{-1}) = -
\infty$. Thus $(0,0)\not\in A^b(w)$ for any $w\in X$ and this refutes (3).

\par If $w_0\in X$, then $A^b(w_0)=X$ because  $\psi = -\infty$ holds
only on $X$. The closure of $X$ contains $(0,0)$ and this means that
$A^b(w_0)$ is not closed and is a proper subset of $A^b((0,0))$ while
$(0,0)\in A^{cb}(w_0)$. This refutes (1) and (2).

\section{Open problems}
\par \be\item Is $A^{cb}(w_0)=A^{cb}_e(w_0)$?
\item Does $A^{cb}(w_0)\sbs \mathbf{c}(M)$?
\item Is there a non-negative $u\in PSH^{cb}(M)$ such that the set
    $\{u=0\}=\mathbf{c}(M)$?
\ee

%
%
 \vspace{1truecm}


\begin{thebibliography}{\;\;\;\;\;}

\bibitem{AS}  A. Aytuna and A. Sadullaev, {\em Parabolic Stein manifolds,}
    Math. Scand. {\bf 114} (2014),  86–-109.

\bibitem{HST1}  T. Harz, N. Shcherbina, and G. Tomassini, {\em On defining
    functions and cores for unbounded domains I,} Math. Z. {\bf 286}
    (2017), 987--1002.

\bibitem{HST2}  T. Harz, N. Shcherbina, and G. Tomassini, {\em On defining
    functions and cores for unbounded domains II,} J. Geom. Anal. (2017),
    First Online: 21 April 2017.

\bibitem{HST3}  T. Harz, N. Shcherbina, and G. Tomassini, {\em On defining
    functions and cores for unbounded domains III,} submitted for
    publication.

\bibitem{HST4}  T. Harz, N. Shcherbina, and G. Tomassini, {\em  Wermer
    type
    sets and extension of CR functions,} Indiana Univ. Math. J. {\bf 61}
    (2012), 431–-459.

\bibitem{NN} R. Narasimhan and Y. Nievergelt, {\it Complex analysis in one
    variable,} 2nd edition. Birkh\"auser, Boston, 2001.

\bibitem{P} E.A. Poletsky, {\em Plurisubharmonic functions as solutions of
    variational problems,} Several complex variables and complex geometry,
    Part 1 (Santa Cruz, CA, 1989), 163–-171, Proc. Sympos. Pure Math., 52,
    Part 1, Amer. Math. Soc., Providence, RI, 1991.

\bibitem{R1} J.P. Rosay, {\em Poletsky theory of disks on holomorphic
    manifolds,} Indiana Univ. Math. J. {\bf 52} (2003), 157–-170.

\bibitem{R2} J.-P. Rosay, {\em Discs in complex manifolds with no bounded
    plurisubharmonic functions,} Proc. Amer. Math. Soc. {\bf 132} (2004),
    2315–-2319.

\bibitem{S1} Z. Slodkowski, {\em Local maximum property and
$q$-plurisubharmonic functions in uniform algebras,}  J. Math. Anal. Appl.
{\bf 115} (1986), 105–-130.

\bibitem{S2} Z. Slodkowski, {\em Preliminary notes,} Spring 2017.

\bibitem{ST} Z. Slodkowski and G. Tomassini, {\em Minimal kernels of
    weakly
    complete spaces,}  J. Funct. Anal. {\bf 210} (2004), 125–-147.


\end{thebibliography}
\end{document}